\documentclass[twoside,11pt]{article}

\usepackage{jmlr2e}
\usepackage{bm}

\usepackage[english]{babel}

\usepackage{dsfont}

\usepackage{amsmath}
\usepackage{arydshln}

\usepackage{amsmath}
\usepackage{blkarray}

\usepackage[final]{pdfpages}

\usepackage[noframe]{showframe}
\usepackage{framed}
\usepackage{lipsum}
\definecolor{color0}  {RGB}{174,225,254} 
\definecolor{color1}  {RGB}{220,227,248} 
\definecolor{color2}  {RGB}{28,130,185} 
\definecolor{color3}  {RGB}{255,253,250} 
{\endMakeFramed}
\definecolor{shadecolor}{gray}{0.75}

\usepackage[framemethod=tikz]{mdframed}
\newcommand{\mdframecolor}{color1}
\newcommand{\mdframehideline}{true}

\usepackage{xcolor}

\usepackage{tcolorbox}

\usepackage{graphicx}  
\usepackage{float}  
\usepackage{subfigure}  

\usepackage{algorithm}
\usepackage{algpseudocode}
\usepackage{amsmath}
\usepackage{graphics}
\usepackage{epsfig}

\usepackage{blkarray}

\usepackage{listings}

\usepackage[colorlinks]{hyperref}

\usepackage{color}   
\usepackage{hyperref}
\hypersetup{
    colorlinks=true, 
    linktoc=all,     
    linkcolor=red,  
}

\newcommand{\real}{\mathbb{R}}

\newcommand{\cspace}{\mathcal{C}}
\newcommand{\nspace}{\mathcal{N}}

\newcommand{\ba}{\bm{a}}
\newcommand{\bA}{\bm{A}}
\newcommand{\bb}{\bm{b}}
\newcommand{\bB}{\bm{B}}
\newcommand{\bc}{\bm{c}}
\newcommand{\bC}{\bm{C}}

\newcommand{\bD}{\bm{D}}

\newcommand{\bE}{\bm{E}}

\newcommand{\bF}{\bm{F}}

\newcommand{\bI}{\bm{I}}

\newcommand{\bL}{\bm{L}}

\newcommand{\bn}{\bm{n}}
\newcommand{\bN}{\bm{N}}

\newcommand{\bP}{\bm{P}}
\newcommand{\bq}{\bm{q}}
\newcommand{\bQ}{\bm{Q}}
\newcommand{\br}{\bm{r}}
\newcommand{\bR}{\bm{R}}

\newcommand{\bT}{\bm{T}}
\newcommand{\bu}{\bm{u}}
\newcommand{\bU}{\bm{U}}
\newcommand{\bv}{\bm{v}}
\newcommand{\bV}{\bm{V}}

\newcommand{\bx}{\bm{x}}

\newcommand{\bZ}{\bm{Z}}

\newcommand{\bzero}{\mathbf{0}}

\jmlrheading{1}{2021}{1-48}{4/00}{10/00}{Jun Lu}

\ShortHeadings{On the Column and Row Ranks of a Matrix}{Jun Lu}
\firstpageno{1}

\begin{document}

\title{On the Column and Row Ranks of a Matrix}

\author{
\begin{center}
\name Jun Lu \\ 
\email jun.lu.locky@gmail.com
\end{center}
   }


\maketitle

\begin{abstract}
Every $m$ by $n$ matrix $\bA$ with rank $r$ has exactly $r$ independent rows and $r$ independent columns. The fact has become the most fundamental theorem in linear algebra such that we may favor it in an unconscious way. The sole aim of this paper is to give a self-contained introduction to concepts and mathematical tools for the rank of a matrix in order to seamlessly introduce how it works in applied linear algebra. However, we clearly realize our inability to cover all the useful and interesting results concerning this topic and given the paucity of scope to present this discussion, e.g., a proof via the injective linear map. We refer the reader to literature in the field of linear algebra for a more detailed introduction to the related fields.


\end{abstract}



\section{The Theorem}\label{append:row-equal-column}

In normal conditions, what we talk about the rank of a matrix $\bA\in \real^{m\times n}$ refers to the dimension of the column space of $\bA$.
The rank of $\bA$ is equal to the maximal number of linearly independent columns of $\bA$. By the following theorem, this is also the maximal number of linearly independent rows of $\bA$. The matrix $\bA$ and its transpose $\bA^\top$ have the same rank. We say that $\bA$ has full rank, if its rank is equal to $min\{m,n\}$. In another word, this is true if and only if either all the columns of $\bA$ are linearly independent, or all the rows of $\bA$ are linearly independent. Specifically, given a vector $\bu \in \real^m$ and a vector $\bv \in \real^n$, then the $m\times n$ matrix $\bu\bv^\top$ is of rank 1.
We here present the theorem rigorously as follows:
\begin{mdframed}[hidealllines=\mdframehideline,backgroundcolor=\mdframecolor]
	\begin{theorem}[Row Rank Equals to Column Rank]\label{lemma:equal-dimension-rank}
		The dimension of the column space of a matrix $\bA\in \real^{m\times n}$ is equal to the dimension of its
		row space, i.e., the row rank and the column rank of a matrix $\bA$ are equal.
	\end{theorem}
\end{mdframed}

This paper offers three ways to prove the equivalence of the row and column ranks:
\begin{itemize}
\item An elementary way that constructs column basis from row basis;

\item A Gram-Schmidt way that results in the UTV decomposition framework which in turn gives insights on the theorem;

\item A column-row decomposition way that arises from Gaussian elimination.
\end{itemize}

\section{An Elementary Way}
In this section, we prove Theorem~\ref{lemma:equal-dimension-rank} by an elementary approach from which it also reveals the fundamental theorem of linear algebra.

\begin{proof}[of Theorem~\ref{lemma:equal-dimension-rank}, \textbf{The First Way: An Elementary Way}]
	We first notice that the null space of $\bA$ \footnote{Null space of $\bA$, denoted as $\nspace(\bA)$; and column space of $\bA$, denoted as $\cspace(\bA)$.} is orthogonal complementary to the row space of $\bA$ ($\cspace(\bA^\top)$, i.e., the column space of $\bA^\top$, where the row space of $\bA$ is exactly the column space of $\bA^\top$): $\nspace(\bA) \bot \cspace(\bA^\top)$, that is, vectors in the null space of $\bA$ are orthogonal to vectors in the row space of $\bA$. To see this, suppose $\bA$ has rows $\ba_1^\top, \ba_2^\top, \ldots, \ba_m^\top$ and $\bA=[\ba_1^\top; \ba_2^\top; \ldots; \ba_m^\top]$ is the row partition of $\bA$. For any vector $\bx\in \nspace(\bA)$, we have $\bA\bx = \bzero$, that is, $[\ba_1^\top\bx; \ba_2^\top\bx; \ldots; \ba_m^\top\bx]=\bzero$. And since the row space of $\bA$ is spanned by $\ba_1^\top, \ba_2^\top, \ldots, \ba_m^\top$. Then $\bx$ is perpendicular to any vectors from $\cspace(\bA^\top)$ which means $\nspace(\bA) \bot \cspace(\bA^\top)$.
	
	Now suppose further that the dimension of row space of $\bA$ is $r$. \textcolor{blue}{Let $\br_1, \br_2, \ldots, \br_r$ be a set of vectors in $\real^n$ and form a basis for the row space of $\bA$}. Then the $r$ vectors $\bA\br_1, \bA\br_2, \ldots, \bA\br_r$ are apparently in the column space of $\bA$,  which are linearly independent. To show this, suppose we have a linear combination of the $r$ vectors: $x_1\bA\br_1 + x_2\bA\br_2+ \ldots+ x_r\bA\br_r=\bzero$, that is, $\bA(x_1\br_1 + x_2\br_2+ \ldots+ x_r\br_r)=\bzero$ and the vector $\bv=x_1\br_1 + x_2\br_2+ \ldots+ x_r\br_r$ is in null space of $\bA$: $\bv\in \nspace(\bA)$. But since $\{\br_1, \br_2, \ldots, \br_r\}$ is a basis for the row space of $\bA$, $\bv$ is thus also in the row space of $\bA$ : $\bv\in \cspace(\bA^\top)$. We have shown above that vectors from null space of $\bA$ is perpendicular to vectors from row space of $\bA$, thus $\bv^\top\bv=0$ and $x_1=x_2=\ldots=x_r=0$. Then \textcolor{blue}{$\bA\br_1, \bA\br_2, \ldots, \bA\br_r$ are in the column space of $\bA$ and they are linearly independent} which means the dimension of the column space of $\bA$ is larger than $r$. This result shows that \textbf{row rank of $\bA\leq $ column rank of $\bA$}. 
	
	If we apply this process again for $\bA^\top$. We will have \textbf{column rank of $\bA\leq $ row rank of $\bA$}. This completes the proof. 
\end{proof}

An important byproduct from the above proof is that if we find a row basis $\{\br_1, \br_2, \ldots, \br_r\}$ of matrix $\bA$ with rank $r$, then, $\{\bA\br_1, \bA\br_2, \ldots, \bA\br_r\}$ is a column basis for it. The reverse is not true normally since $\bA$ may not be invertible. 

However, suppose $\{\bc_1, \bc_2, \ldots, \bc_r\}$ is a column basis for $\bA$, i.e., $\{\bc_1, \bc_2, \ldots, \bc_r\}$ is a row basis for $\bA^\top$. Apply the finding above, it is trivial that $\{\bA^\top\bc_1, \bA^\top\bc_2, \ldots, \bA^\top\bc_r\}$ is a row basis for $\bA$ (i.e., a column basis for $\bA^\top$). 

Specially, when $\bA\in \real^{n\times n}$ is nonsingular and $\{\bc_1, \bc_2, \ldots, \bc_n\}$ is a column basis for $\bA$, then $\{\bA^{-1}\bc_1, \bA^{-1}\bc_2, $ $\ldots, $ $\bA^{-1}\bc_n\}$ is a row basis for it. Or apparently, in this case, the $n$ rows or even the $n$ columns of $\bA$ also form a row basis as well since both the rows and columns span the whole $\real^n$ space .

\subsection{Fundamental Theorem of Linear Algebra}\label{appendix:fundamental-rank-nullity}
The byproduct of the above proof is essential to prove the fundamental theorem of linear algebra. As a recap, for any matrix $\bA\in \real^{m\times n}$, it can be easily verified that any vector in row space of $\bA$ is perpendicular to any vector in null space of $\bA$. Suppose $\bx_n \in \nspace(\bA)$, then $\bA\bx_n = \bzero$ such that $\bx_n$ is perpendicular to every row of $\bA$ which agrees with our claim. 

Similarly, we can also show that any vector in column space of $\bA$ is perpendicular to any vector in null space of $\bA^\top$. Further, the column space of $\bA$ and the null space of $\bA^\top$ span the whole $\real^m$ space which is known as the fundamental theorem of linear algebra.

The fundamental theorem contains two parts, the dimension of the subspaces and the orthogonality of the subspaces. The orthogonality can be easily verified as shown above. Moreover, when the row space has dimension $r$, the null space has dimension $n-r$. This cannot be easily stated and we prove in the following theorem.

\begin{figure}[h!]
	\centering
	\includegraphics[width=0.98\textwidth]{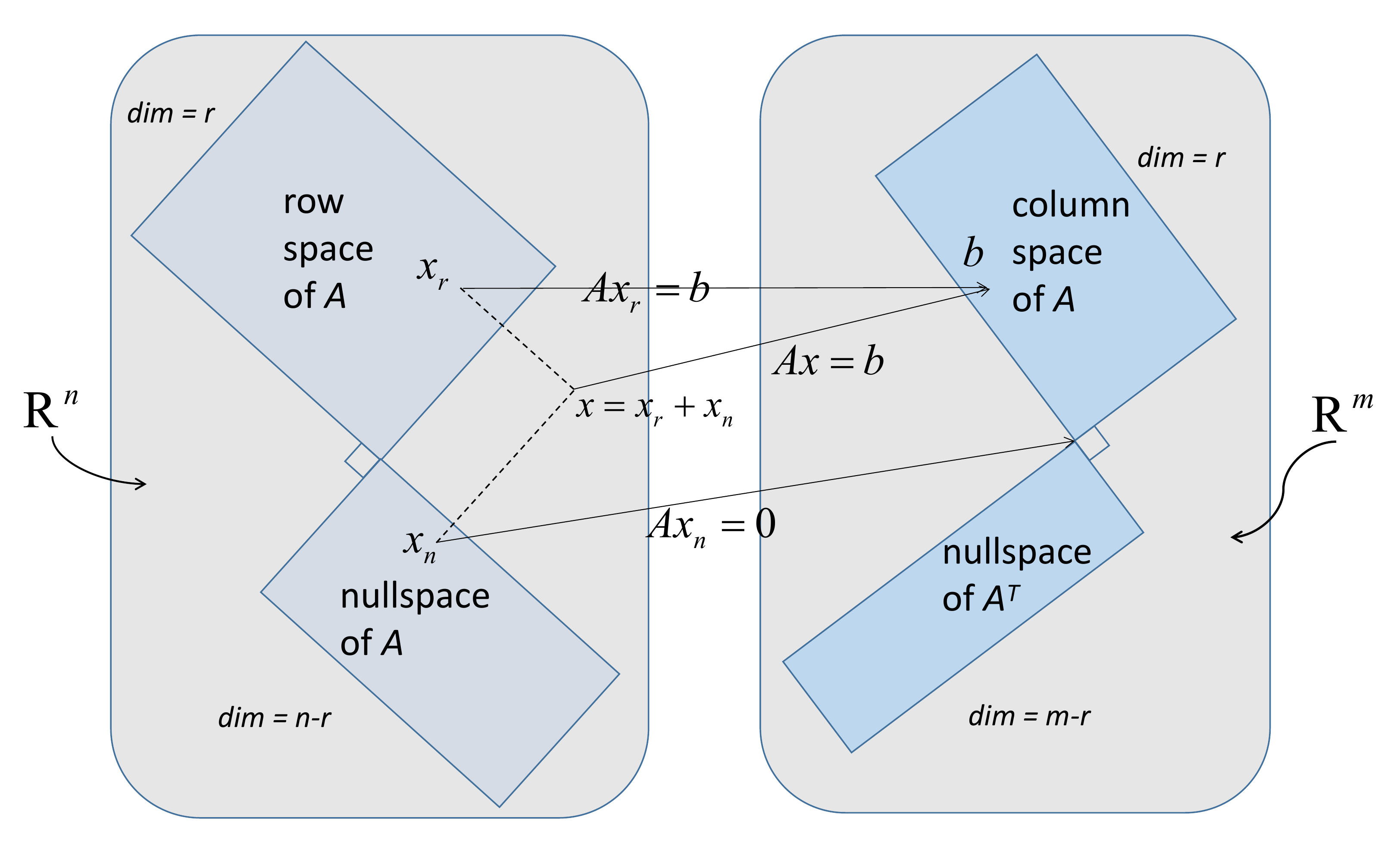}
	\caption{Two pairs of orthogonal subspaces in $\real^n$ and $\real^m$. $\dim(\cspace(\bA^\top)) + \dim(\nspace(\bA))=n$ and $\dim(\nspace(\bA^\top)) + \dim(\cspace(\bA))=m$. The null space component goes to zero as $\bA\bx_n = \bzero \in \real^m$. The row space component goes to column space as $\bA\bx_r = \bA(\bx_r+\bx_n)=\bb \in \cspace(\bA)$.}
	\label{fig:lafundamental}
\end{figure}
\begin{mdframed}[hidealllines=\mdframehideline,backgroundcolor=\mdframecolor]
	\begin{theorem}[Fundamental Theorem of Linear Algebra]\label{theorem:fundamental-linear-algebra}
		
		Orthogonal Complement and Rank-Nullity Theorem: for any matrix $\bA\in \real^{m\times n}$, we have 
		
		$\bullet$ $\nspace(\bA)$ is orthogonal complement to the row space $\cspace(\bA^\top)$ in $\real^n$: $dim(\nspace(\bA))+dim(\cspace(\bA^\top))=n$;
		
		$\bullet$ $\nspace(\bA^\top)$ is orthogonal complement to the column space $\cspace(\bA)$ in $\real^m$: $dim(\nspace(\bA^\top))+dim(\cspace(\bA))=m$;
		
		$\bullet$ For rank-$r$ matrix $\bA$, $dim(\cspace(\bA^\top)) = dim(\cspace(\bA)) = r$, that is, $dim(\nspace(\bA)) = n-r$ and $dim(\nspace(\bA^\top))=m-r$.
	\end{theorem}
\end{mdframed}

\begin{proof}[of Theorem~\ref{theorem:fundamental-linear-algebra}]
	Follow from the proof of Theorem~\ref{lemma:equal-dimension-rank}. Let $\br_1, \br_2, \ldots, \br_r$ be a set of vectors in $\real^n$ that form a basis for the row space, then $\bA\br_1, \bA\br_2, \ldots, \bA\br_r$ is a basis for the column space of $\bA$. Let $\bn_1, \bn_2, \ldots, \bn_k \in \real^n$ form a basis for the null space of $\bA$. Follow again from the proof of Theorem~\ref{lemma:equal-dimension-rank}, $\nspace(\bA) \bot \cspace(\bA^\top)$, thus, $\br_1, \br_2, \ldots, \br_r$ are perpendicular to $\bn_1, \bn_2, \ldots, \bn_k$. Then, $\{\br_1, \br_2, \ldots, \br_r, \bn_1, \bn_2, \ldots, \bn_k\}$ is linearly independent in $\real^n$.
	
	For any vector $\bx\in \real^n $, $\bA\bx$ is in the column space of $\bA$. Then it can be expressed by a combination of $\bA\br_1, \bA\br_2, \ldots, \bA\br_r$: $\bA\bx = \sum_{i=1}^{r}a_i\bA\br_i$ which states that $\bA(\bx-\sum_{i=1}^{r}a_i\br_i) = \bzero$ and $\bx-\sum_{i=1}^{r}a_i\br_i$ is thus in $\nspace(\bA)$. Since $\{\bn_1, \bn_2, \ldots, \bn_k\}$ is a basis for the null space of $\bA$, $\bx-\sum_{i=1}^{r}a_i\br_i$ can be expressed by a combination of $\bn_1, \bn_2, \ldots, \bn_k$: $\bx-\sum_{i=1}^{r}a_i\br_i = \sum_{j=1}^{k}b_j \bn_j$, i.e., $\bx=\sum_{i=1}^{r}a_i\br_i + \sum_{j=1}^{k}b_j \bn_j$. That is, any vector $\bx\in \real^n$ can be expressed by $\{\br_1, \br_2, \ldots, \br_r, \bn_1, \bn_2, \ldots, \bn_k\}$ and the set forms a basis for $\real^n$. Thus the dimension sum to $n$: $r+k=n$, i.e., $dim(\nspace(\bA))+dim(\cspace(\bA^\top))=n$. Similarly, we can prove $dim(\nspace(\bA^\top))+dim(\cspace(\bA))=m$.
\end{proof}

Figure~\ref{fig:lafundamental} demonstrates two pairs of such orthogonal subspaces and shows how $\bA$ takes $\bx$ into the column space. The dimensions of row space of $\bA$ and null space of $\bA$ add to $n$. And the dimensions of column space of $\bA$ and null space of $\bA^\top$ add to $m$. The null space component goes to zero as $\bA\bx_{\bn} = \bzero \in \real^m$ which is the intersection of column space of $\bA$ and null space of $\bA^\top$. The row space component goes to column space as $\bA\bx_{\br} = \bA(\bx_{\br} + \bx_{\bn})=\bb\in \real^m$. 

The fundamental theorem of linear algebra was discussed in \citep{strang1993fundamental}, and a more detailed review is provided in \citep{lu2021revisit} where the authors provide 7 figures from different perspectives (e.g., elementary perspective, least squares, pseudo-inverse and SVD) to describe the theorem.

\section{Gram-Schmidt Way}\label{section:gaussian-way-main}
The UTV decomposition goes further from the QR decomposition which in tern relies on the \textit{Gram-Schmidt process} by factoring the matrix into two orthogonal matrices $\bA=\bU\bT\bV$, where $\bU, \bV$ are orthogonal. The resulting $\bT$ supports rank estimation. The matrix $\bT$ can be lower triangular which results in the ULV decomposition, or it can be upper triangular which results in the URV decomposition. The UTV framework shares similar form as the singular value decomposition (SVD) and can be regarded as inexpensive alternatives to the SVD.

We here firstly formulate the decomposition in the following theorem, and the proof of the existence of the ULV decomposition will be discussed in Section~\ref{section:utv-prove}. \footnote{Notice that original claim for the full ULV decomposition relies on the rank of matrix $\bA$ to be $r$ rather than the column rank \citep{lu2021numerical}. But here we slightly modify the argument since we want to employ it into the proof of the row rank equals to column.}
\begin{mdframed}[hidealllines=\mdframehideline,backgroundcolor=\mdframecolor]
	\begin{theorem}[Full ULV Decomposition]\label{theorem:ulv-decomposition}
		Every $m\times n$ matrix $\bA$ with \textit{column rank} $r$ can be factored as 
		$$
		\bA = \bU \begin{bmatrix}
			\bL & \bzero \\
			\bzero & \bzero 
		\end{bmatrix}\bV,
		$$
		where $\bU\in \real^{m\times m}$ and $\bV\in \real^{n\times n}$ are two orthogonal matrices, and $\bL\in \real^{r\times r}$ is a lower triangular matrix.
	\end{theorem}
\end{mdframed}

From the ULV decomposition, we here provide an alternative way to prove the column rank and row rank of matrix $\bA$ are the same. Note that, in the proof of the QR, LQ and UTV decompositions, we do not utilize the results that the column and row ranks of a matrix are the same. Otherwise, we won't be able to prove it via the existence of the ULV decomposition. 

\begin{proof}[of Theorem~\ref{lemma:equal-dimension-rank}, \textbf{A Second Way}]
	Any $m\times n$ matrix $\bA$ with column rank $r$ can be factored as 
	$$
	\bA = \bU_0 \begin{bmatrix}
		\bL & \bzero \\
		\bzero & \bzero 
	\end{bmatrix}\bV_0,
	$$
	where $\bU_0\in \real^{m\times m}$ and $\bV_0\in \real^{n\times n}$ are two orthogonal matrices, and $\bL\in \real^{r\times r}$ is a lower triangular matrix. Let $\bD = \begin{bmatrix}
		\bL & \bzero \\
		\bzero & \bzero 
	\end{bmatrix}$, the row rank and column rank of $\bD$ are apparently the same. If we could prove the column rank of $\bA$ equals to the column rank of $\bD$, and the row rank of $\bA$ equals to the row rank of $\bD$, then we complete the proof. 
	
	Let $\bU = \bU_0^\top$, $\bV=\bV_0^\top$, then it follows that $\bD = \bU\bA\bV$. Decompose the above idea into two steps, a moment of reflexion reveals that, 
	\begin{itemize}
		\item if we could first prove the row rank and column rank of $\bA$ are equal to those of $\bU\bA$; 
		\item and then, if we further prove the row rank and column rank of $\bU\bA $ are equal to those of $\bU\bA\bV=\bD$, we could also complete the proof.
	\end{itemize}
	
	\paragraph{Row rank and column ranks of $\bA$ are equal to those of $\bU\bA$} Let $\bB = \bU\bA$, and let further $\bA=[\ba_1,\ba_2,\ldots, \ba_n]$ and $\bB=[\bb_1,\bb_2,\ldots,\bb_n]$ be the column partitions of $\bA$ and $\bB$. Therefore, $[\bb_1,\bb_2,\ldots,\bb_n] = [\bU\ba_1,\bU\ba_2,\ldots, \bU\ba_n]$. If $x_1\ba_1+x_2\ba_2+\ldots+x_n\ba_n=0$, then we also obtain 
	$$
	\bU(x_1\ba_1+x_2\ba_2+\ldots+x_n\ba_n) = x_1\bb_1+x_2\bb_2+\ldots+x_n\bb_n = 0.
	$$
	Let $j_1, j_2, \ldots, j_r$ be distinct indices between 1 and $n$, if the set $\{\ba_{j_1}, \ba_{j_2}, \ldots, \ba_{j_r}\}$ is linearly independent, the set $\{\bb_{j_1}, \bb_{j_2}, \ldots, \bb_{j_r}\}$ must also be independent. This implies 
	$$
	dim(\cspace(\bB)) \leq dim(\cspace(\bA)).
	$$
	Similarly, by $\bA = \bU^\top\bB$, we have 
	$$
	dim(\cspace(\bA)) \leq dim(\cspace(\bB)).
	$$
	This implies 
	$$
	dim(\cspace(\bB)) = dim(\cspace(\bA)).
	$$
	Apply the process onto $\bB^\top$ and $\bA^\top$, we have 
	$$
	dim(\cspace(\bB^\top)) = dim(\cspace(\bA^\top)).
	$$
	This implies the row rank and column rank of $\bA$ and $\bB=\bU\bA$ are the same.
	Similarly, we can also show that the row rank and column rank of $\bU\bA$ and $\bU\bA\bV$ are the same. This completes the proof.
\end{proof}

\paragraph{Warning on the proof} We might ask: since the UTV framework is an inexpansive alternatives to the SVD, why do not we use the SVD to prove what we want? The answer needs closer scrutiny. The existence of the SVD relies on the fact that the $\bA^\top\bA$ and $\bA$ have same rank which in turn is based on the fundamental theorem of linear algebra that is a byproduct of Theorem~\ref{lemma:equal-dimension-rank} as we have shown in last section. Therefore, this is problematic as we use claim $A$ to prove claim $B$ where $A$ relies on $B$ as well.

\subsection{The Proof of the UTV Decomposition}\label{section:utv-prove}
The most important application of the Gram-Schmidt process is to prove the existence of the QR decomposition, which in turn is essential to prove Theorem~\ref{lemma:equal-dimension-rank}. In this section, we first present the result of the QR decomposition and provide the proof of it via the Gram-Schmidt process, of which other proofs also exist that might be more numerical stable, e.g., via the Householder reflectors and Givens rotations. The detailed discussion can be found in \citep{trefethen1997numerical, golub2013matrix,van2020advanced, lu2021numerical}.  

\subsubsection{QR Decomposition}
In many applications, we are interested in the column space of a matrix $\bA=[\ba_1, \ba_2, ..., \ba_n] \in \real^{m\times n}$. The successive spaces spanned by the columns $\ba_1, \ba_2, \ldots$ of $\bA$ are
$$
\cspace([\ba_1])\,\,\,\, \subseteq\,\,\,\, \cspace([\ba_1, \ba_2]) \,\,\,\,\subseteq\,\,\,\, \cspace([\ba_1, \ba_2, \ba_3])\,\,\,\, \subseteq\,\,\,\, \ldots,
$$
where $\cspace([\ldots])$ is the subspace spanned by the vectors included in the brackets. The idea of the QR decomposition is the construction of a sequence of orthonormal vectors $\bq_1, \bq_2, \ldots$ that span the same successive subspaces. 
$$
\cspace([\bq_1])=\cspace([\ba_1])\,\,\,\, \subseteq\,\,\,\, \cspace([\bq_1, \bq_2])=\cspace([\ba_1, \ba_2]) \,\,\,\,\subseteq\,\,\,\, \cspace([\bq_1, \bq_2, \bq_3])=\cspace([\ba_1, \ba_2, \ba_3])\,\,\,\, \subseteq\,\,\,\, \ldots,
$$
We provide the result of QR decomposition in the following theorem and we delay the discussion of the existence of it in the next sections.

\begin{mdframed}[hidealllines=\mdframehideline,backgroundcolor=\mdframecolor]
	\begin{theorem}[QR Decomposition]
		Every $m\times n$ matrix $\bA=[\ba_1, \ba_2, ..., \ba_n]$ with $m\geq n$ can be factored as 
		$$
		\bA = \bQ\bR,
		$$
		where 
		
		1. \textbf{Reduced}: $\bQ$ is $m\times n$ with orthonormal columns and $\bR$ is an $n\times n$ upper triangular matrix, which is known as the \textbf{reduced QR decomposition};
		
		2. \textbf{Full}: $\bQ$ is $m\times m$ with orthonormal columns and $\bR$ is an $m\times n$ upper triangular matrix which is known as the \textbf{full QR decomposition}. If further restrict the upper triangular matrix to be a square matrix, the full QR decomposition can be denoted as 
		$$
		\bA = \bQ\begin{bmatrix}
			\bR_0\\
			\bzero
		\end{bmatrix},
		$$
		where $\bR_0$ is an $m\times m$ upper triangular matrix. 
		
		Specifically, when $\bA$ has full rank, i.e., has independent columns, $\bR$ also has independent columns, and $\bR$ is nonsingular for the reduced case (i.e., with nonzero diagonals).
	\end{theorem}
\end{mdframed}
The QR decomposition can be extended to matrices having dependent columns. In this case, the diagonals of $\bR$ might contain zeros. But for the proof of the UTV decomposition and the Theorem~\ref{append:row-equal-column}, we only discuss the matrix with independent columns (i.e., full column rank).

\subsubsection{Project a Vector Onto Another Vector}\label{section:project-onto-a-vector}
\begin{figure}[h!]
	\centering  
	\vspace{-0.35cm} 
	\subfigtopskip=2pt 
	\subfigbottomskip=2pt 
	\subfigcapskip=-5pt 
	\subfigure[Project onto a line]{\label{fig:project-line}
		\includegraphics[width=0.47\linewidth]{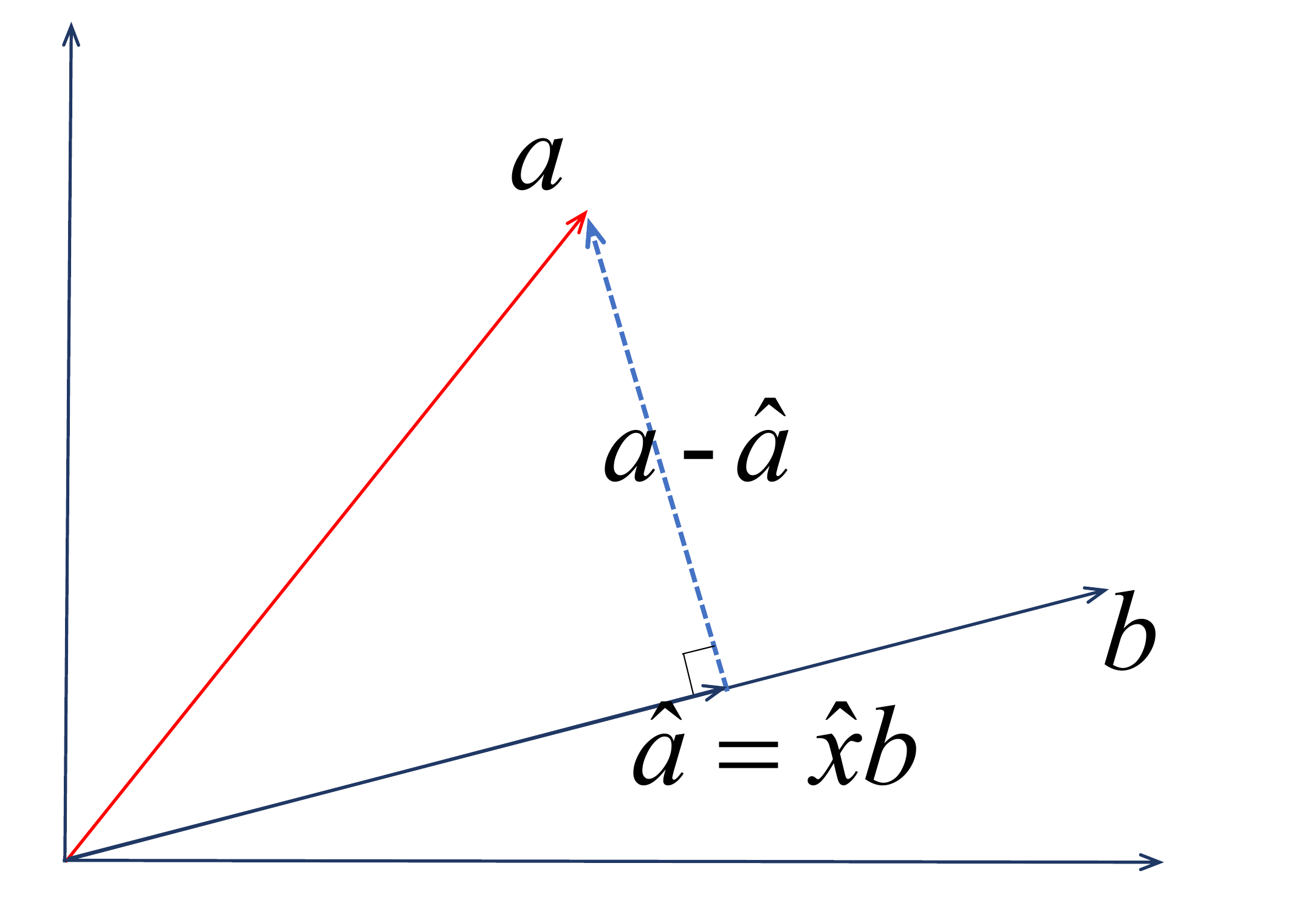}}
	\quad 
	\subfigure[Project onto a space]{\label{fig:project-space}
		\includegraphics[width=0.47\linewidth]{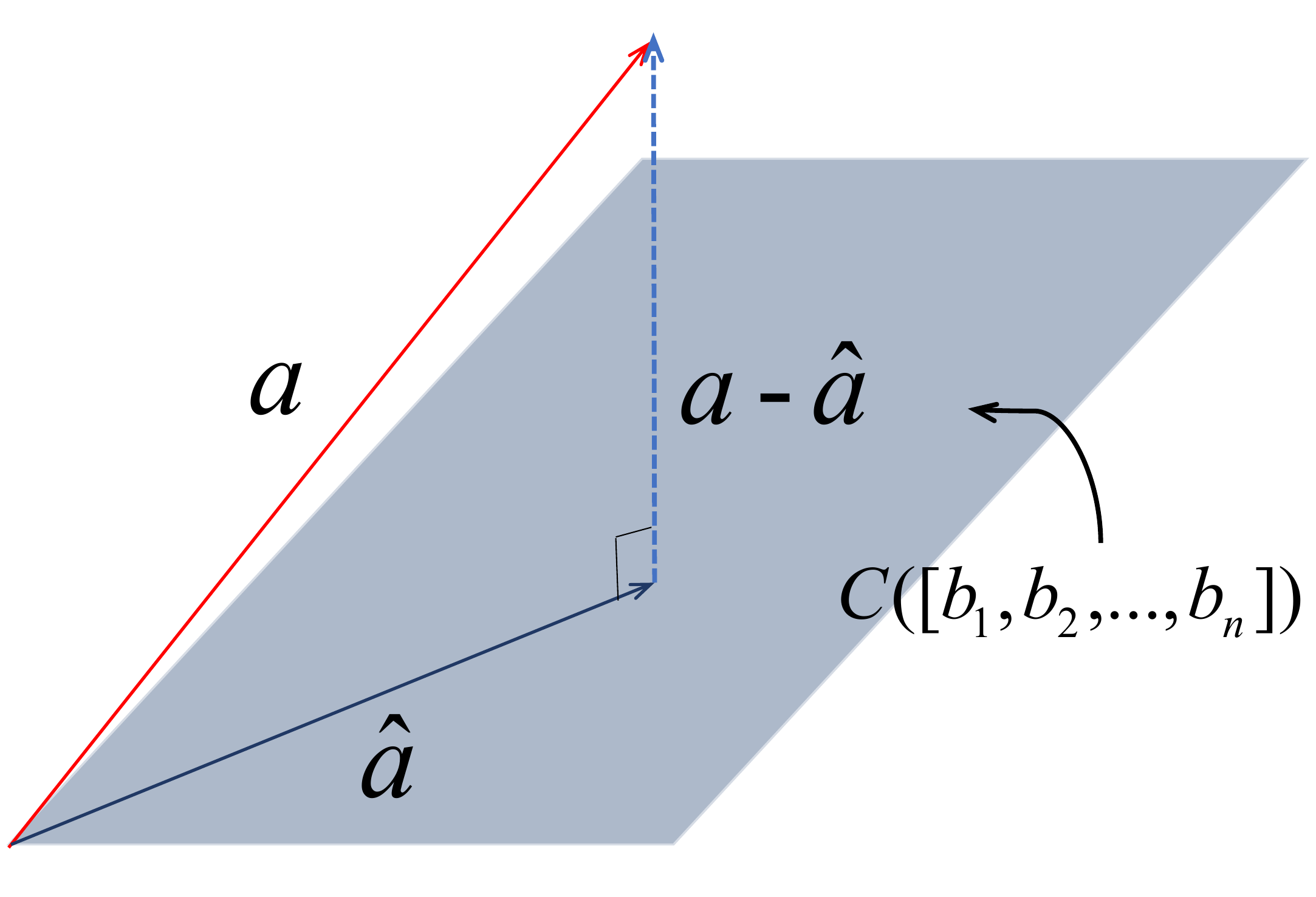}}
	\caption{Project a vector onto a line and onto a space.}
	\label{fig:projection-qr}
\end{figure}

Project a vector $\ba$ to a vector $\bb$ is to find the vector closest to $\ba$ on the line of $\bb$. The projection vector $\hat{\ba}$ is some multiple of $\bb$. Let $\hat{\ba} = \hat{x} \bb$ and $\ba-\hat{\ba}$ is perpendicular to $\bb$ as shown in Figure~\ref{fig:project-line}. We then get the following result:
\begin{mdframed}[hidealllines=\mdframehideline,backgroundcolor=\mdframecolor,frametitle={Project Vector $\ba$ Onto Vector $\bb$}]
	$\ba-\hat{\ba}$ is perpendicular to $\bb$, so $(\ba-\hat{x}\bb)^\top\bb=0$: $\hat{x}$ = $\frac{\ba^\top\bb}{\bb^\top\bb}$ and $\hat{\ba} = \frac{\ba^\top\bb}{\bb^\top\bb}\bb = \frac{\bb\bb^\top}{\bb^\top\bb}\ba$.
\end{mdframed}
%

Further, projecting a vector $\ba$ to a space spanned by $\bb_1, \bb_2, \ldots, \bb_n$ is to find the vector closest to $\ba$ on the column space of $[\bb_1, \bb_2, \ldots, \bb_n]$ can be cast in a similar way. The projection vector $\hat{\ba}$ is a combination of $\bb_1, \bb_2, \ldots, \bb_n$: $\hat{\ba} = \hat{x}_1\bb_1+ \hat{x}_2\bb_2+\ldots+\hat{x}_n\bb_n$. This is actually a least squares problem. To find the projection, we just solve the normal equation $\bB^\top\bB\hat{\bx} = \bB^\top\ba$ where $\bB=[\bb_1, \bb_2, \ldots, \bb_n]$ and $\hat{\bx}=[\hat{x}_1, \hat{x}_2, \ldots, \hat{x}_n]$. We refer the details of this section to \citep{strang1993introduction, trefethen1997numerical, golub2013matrix, lu2021rigorous} as it is not the main interest of this paper.

\subsubsection{Existence of the QR Decomposition via the Gram-Schmidt Process}\label{section:gram-schmidt-process}
For three independent vectors $\ba_1, \ba_2, \ba_3$ and the space spanned by the three independent vectors $\cspace{([\ba_1, \ba_2, \ba_3])}$, i.e., the column space of matrix $[\ba_1, \ba_2, \ba_3]$. We intend to construct three orthogonal vectors $\bb_1, \bb_2, \bb_3$ in which case $\cspace{([\bb_1, \bb_2, \bb_3])}$ = $\cspace{([\ba_1, \ba_2, \ba_3])}$. Then we divide the orthogonal vectors by their length to normalize. This process produces three orthonormal vectors $\bq_1 = \frac{\bb_1}{||\bb_1||}$, $\bq_2 = \frac{\bb_2}{||\bb_2||}$, $\bq_2 = \frac{\bb_2}{||\bb_2||}$.

\begin{figure}[H]
	\centering  
	\vspace{-0.35cm} 
	\subfigtopskip=2pt 
	\subfigbottomskip=2pt 
	\subfigcapskip=-5pt 
	\subfigure[Project $\ba_2$ onto the space perpendicular to $\bb_1$.]{\label{fig:gram-schmidt1}
		\includegraphics[width=0.47\linewidth]{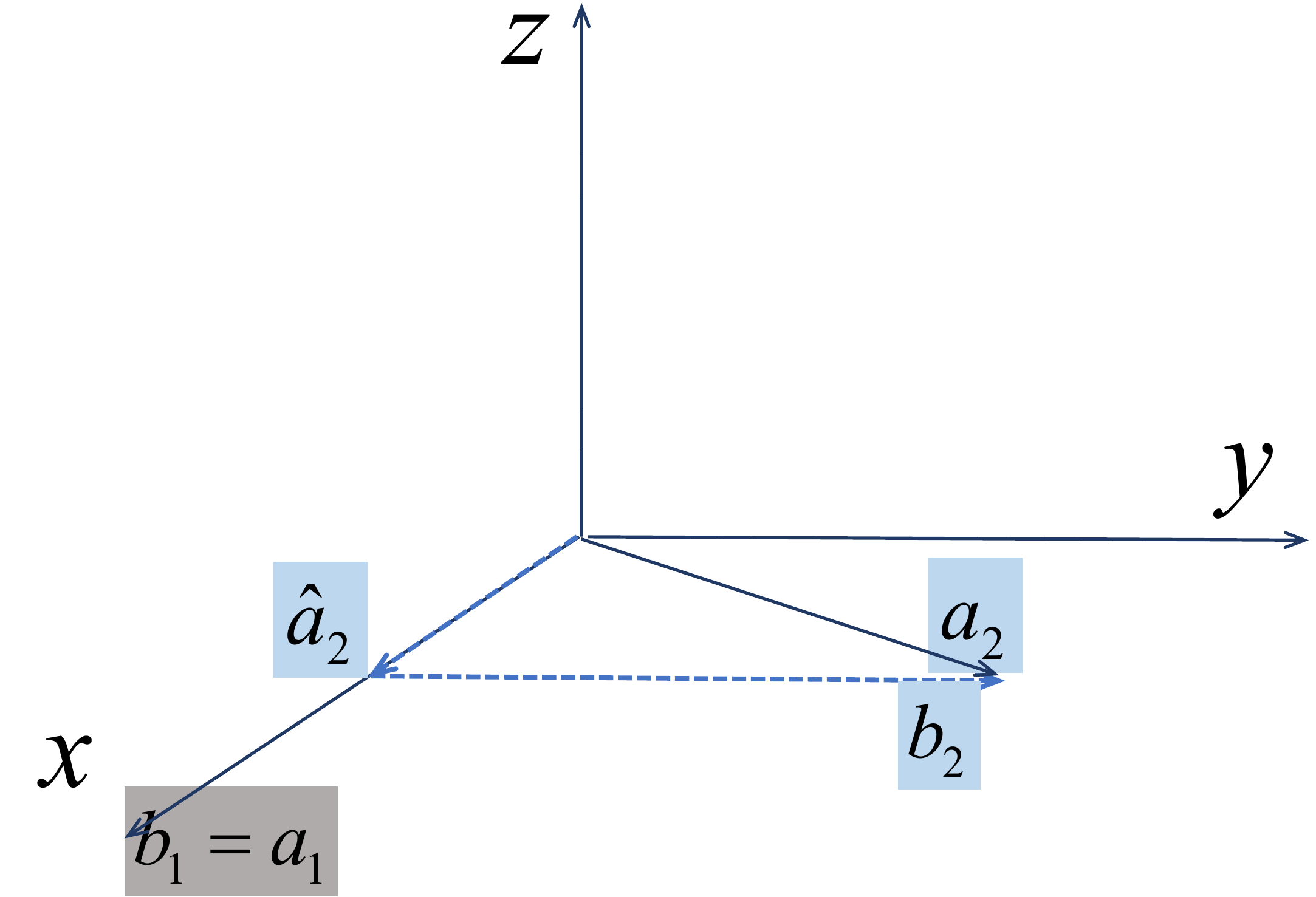}}
	\quad 
	\subfigure[Project $\ba_3$ onto the space perpendicular to $\bb_1, \bb_2$.]{\label{fig:gram-schmidt2}
		\includegraphics[width=0.47\linewidth]{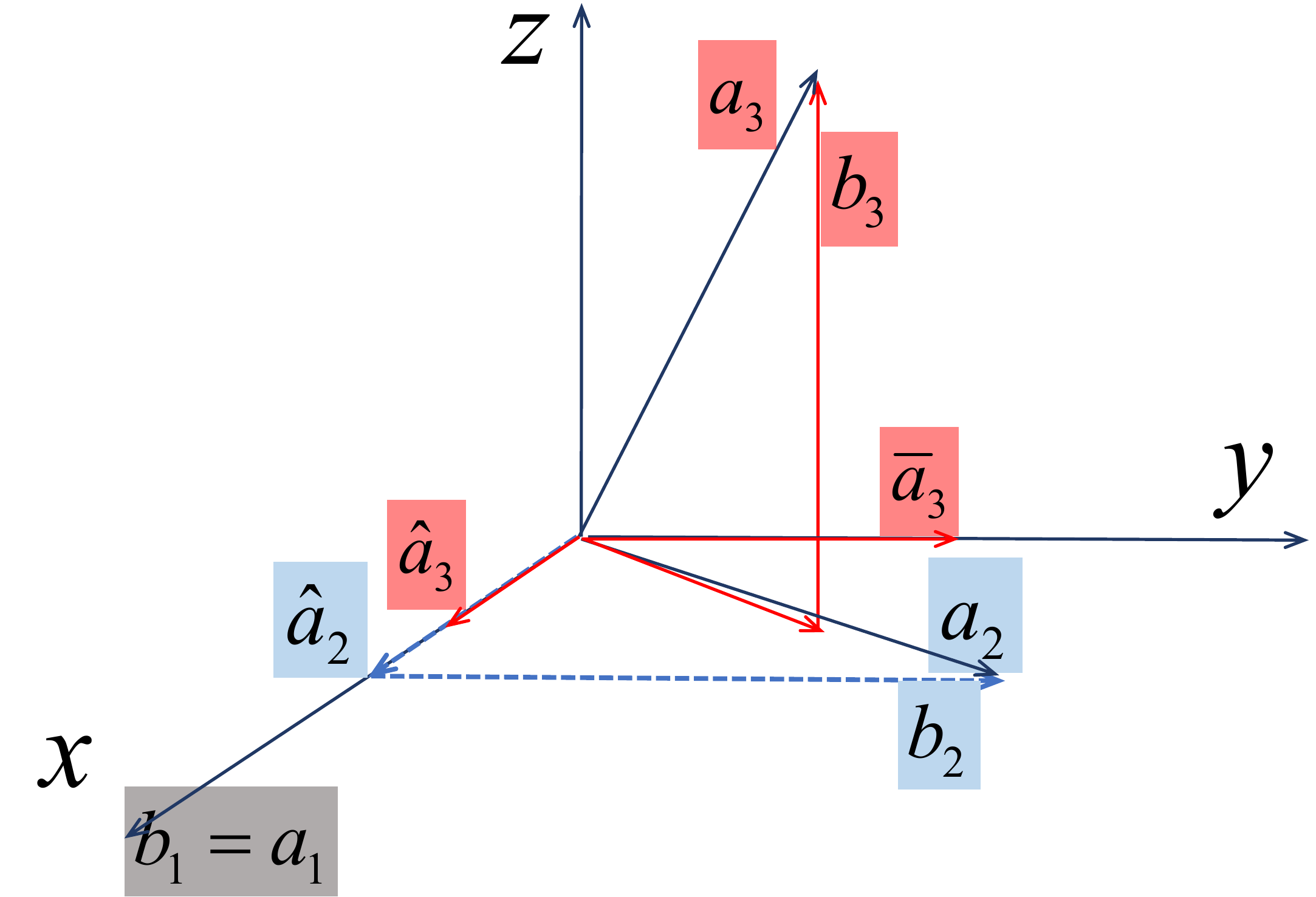}}
	\caption{Gram-Schmidt process.}
	\label{fig:gram-schmidt-12}
\end{figure}

For the first vector, we choose $\bb_1 = \ba_1$ directly. The second vector $\bb_1$ must be perpendicular to the first one. This is actually the vector $\ba_2$ subtracting its projection along $\bb_1$:
\begin{equation}
	\begin{aligned}
		\bb_2 &= \ba_2- \frac{\bb_1 \bb_1^\top}{\bb_1^\top\bb_1} \ba_2 = (\bI- \frac{\bb_1 \bb_1^\top}{\bb_1^\top\bb_1} )\ba_2   \qquad &(\text{Projection view})\\
		&= \ba_2-  \underbrace{\frac{ \bb_1^\top \ba_2}{\bb_1^\top\bb_1} \bb_1}_{\hat{\ba}_2}, \qquad &(\text{Combination view}) \nonumber
	\end{aligned}
\end{equation}
where the first equation shows $\bb_2$ is a multiplication of a matrix $(\bI- \frac{\bb_1 \bb_1^\top}{\bb_1^\top\bb_1} )$ and $\ba_2$, i.e., project $\ba_2$ onto the orthogonal complement space of $\cspace{([\bb_1])}$. The second equation shows $\ba_2$ is a combination of $\bb_1$ and $\bb_2$.
Clearly, the space spanned by $\bb_1, \bb_2$ is the same space spanned by $\ba_1, \ba_2$. The siguation is shown in Figure~\ref{fig:gram-schmidt1} in which we choose \textit{the direction of $\bb_1$ as the $x$-axis in the Cartesian coordinate system}. $\hat{\ba}_2$ is the projection of $\ba_2$ onto the line $\bb_1$. It can be clearly shown that the part of $\ba_2$ perpendicular to $\bb_1$ is $\bb_2 = \ba_2 - \hat{\ba}_2$ from the figure.

For the third vector $\bb_3$, it must be perpendicular to both the $\bb_1$ and $\bb_2$ which is actually the vector $\ba_3$ subtracting its projection along the plane spanned by $\bb_1$ and $\bb_2$
\begin{equation}\label{equation:gram-schdt-eq2}
	\begin{aligned}
		\bb_3 &= \ba_3- \frac{\bb_1 \bb_1^\top}{\bb_1^\top\bb_1} \ba_3 - \frac{\bb_2 \bb_2^\top}{\bb_2^\top\bb_2} \ba_3 = (\bI- \frac{\bb_1 \bb_1^\top}{\bb_1^\top\bb_1}  - \frac{\bb_2 \bb_2^\top}{\bb_2^\top\bb_2} )\ba_3   \qquad &(\text{Projection view})\\
		&= \ba_3- \underbrace{\frac{ \bb_1^\top\ba_3}{\bb_1^\top\bb_1} \bb_1}_{\hat{\ba}_3} - \underbrace{\frac{ \bb_2^\top\ba_3}{\bb_2^\top\bb_2}  \bb_2}_{\bar{\ba}_3},    \qquad &(\text{Combination view})
	\end{aligned}
\end{equation}
where the first equation shows $\bb_3$ is a multiplication of a matrix $(\bI- \frac{\bb_1 \bb_1^\top}{\bb_1^\top\bb_1}  - \frac{\bb_2 \bb_2^\top}{\bb_2^\top\bb_2} )$ and $\ba_3$, i.e., project $\ba_3$ onto the orthogonal complement space of $\cspace{([\bb_1, \bb_2])}$. The second equation shows $\ba_3$ is a combination of $\{\bb_1, \bb_2, \bb_3\}$. We will see this property is essential in the idea of the QR decomposition.
Again, it can be shown that the space spanned by $\{\bb_1, \bb_2, \bb_3\}$ is the same space spanned by $\{\ba_1, \ba_2, \ba_3\}$. The situation is shown in Figure~\ref{fig:gram-schmidt2}, in which we choose \textit{the direction of $\bb_2$ as the $y$-axis of the Cartesian coordinate system}. $\hat{\ba}_3$ is the projection of $\ba_3$ onto the line $\bb_1$, $\bar{\ba}_3$ is the projection of $\ba_3$ onto the line $\bb_2$. It can be shown that the part of $\ba_3$ perpendicular to both $\bb_1$ and $\bb_2$ is $\bb_3=\ba_3-\hat{\ba}_3-\bar{\ba}_3$ from the figure.

Finally, we normalize each vector by dividing their length which produces three orthonormal vectors $\bq_1 = \frac{\bb_1}{||\bb_1||}$, $\bq_2 = \frac{\bb_2}{||\bb_2||}$, $\bq_2 = \frac{\bb_2}{||\bb_2||}$.

This idea can be extended to a set of vectors rather than only three. And we call this process as \textit{Gram-Schmidt process}. After this process, matrix $\bA$ will be triangularized. 

As we mentioned previously, the idea of the QR decomposition is the construction of a sequence of orthonormal vectors $\bq_1, \bq_2, \ldots$ that span the same successive subspaces. 
$$
\cspace([\bq_1])=\cspace([\ba_1])\,\,\,\, \subseteq\,\,\,\, \cspace([\bq_1, \bq_2])=\cspace([\ba_1, \ba_2]) \,\,\,\,\subseteq\,\,\,\, \cspace([\bq_1, \bq_2, \bq_3])=\cspace([\ba_1, \ba_2, \ba_3])\,\,\,\, \subseteq\,\,\,\, \ldots,
$$
This implies any $\ba_k$ is in the space spanned by $\cspace([\bq_1, \bq_2, \ldots, \bq_k])$. \footnote{And also, any $\bq_k$ is in the space spanned by $\cspace([\ba_1, \ba_2, \ldots, \ba_k])$.} As long as we have found these orthonormal vectors, to reconstruct $\ba_i$'s from the orthonormal matrix $\bQ=[\bq_1, \bq_2, \ldots, \bq_n]$, an upper triangular matrix $\bR$ is needed such that $\bA = \bQ\bR$.

\begin{figure}[H]
	\centering  
	\vspace{-0.35cm} 
	\subfigtopskip=2pt 
	\subfigbottomskip=2pt 
	\subfigcapskip=-5pt 
	\subfigure[Reduced QR decomposition]{\label{fig:gphalf}
		\includegraphics[width=0.47\linewidth]{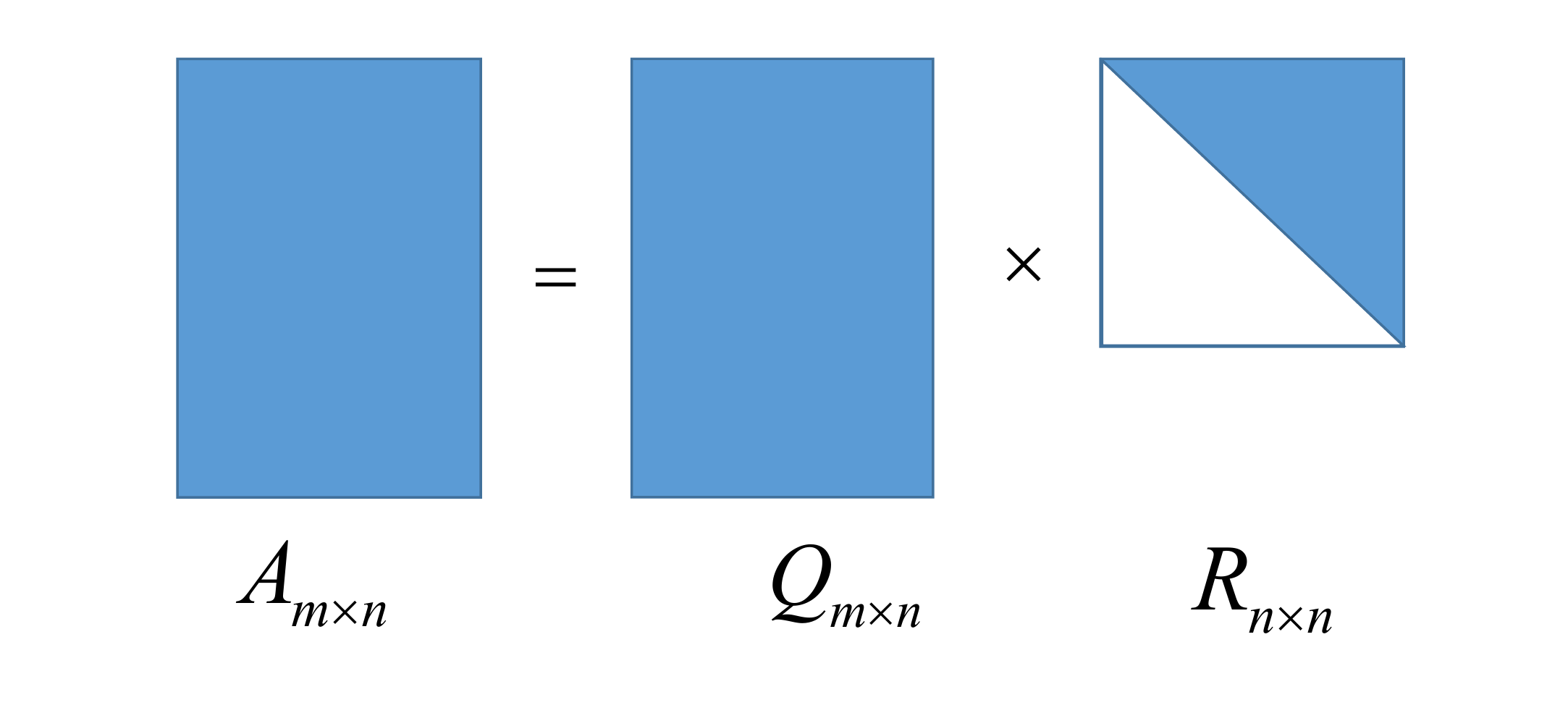}}
	\quad 
	\subfigure[Full QR decomposition]{\label{fig:gpall}
		\includegraphics[width=0.47\linewidth]{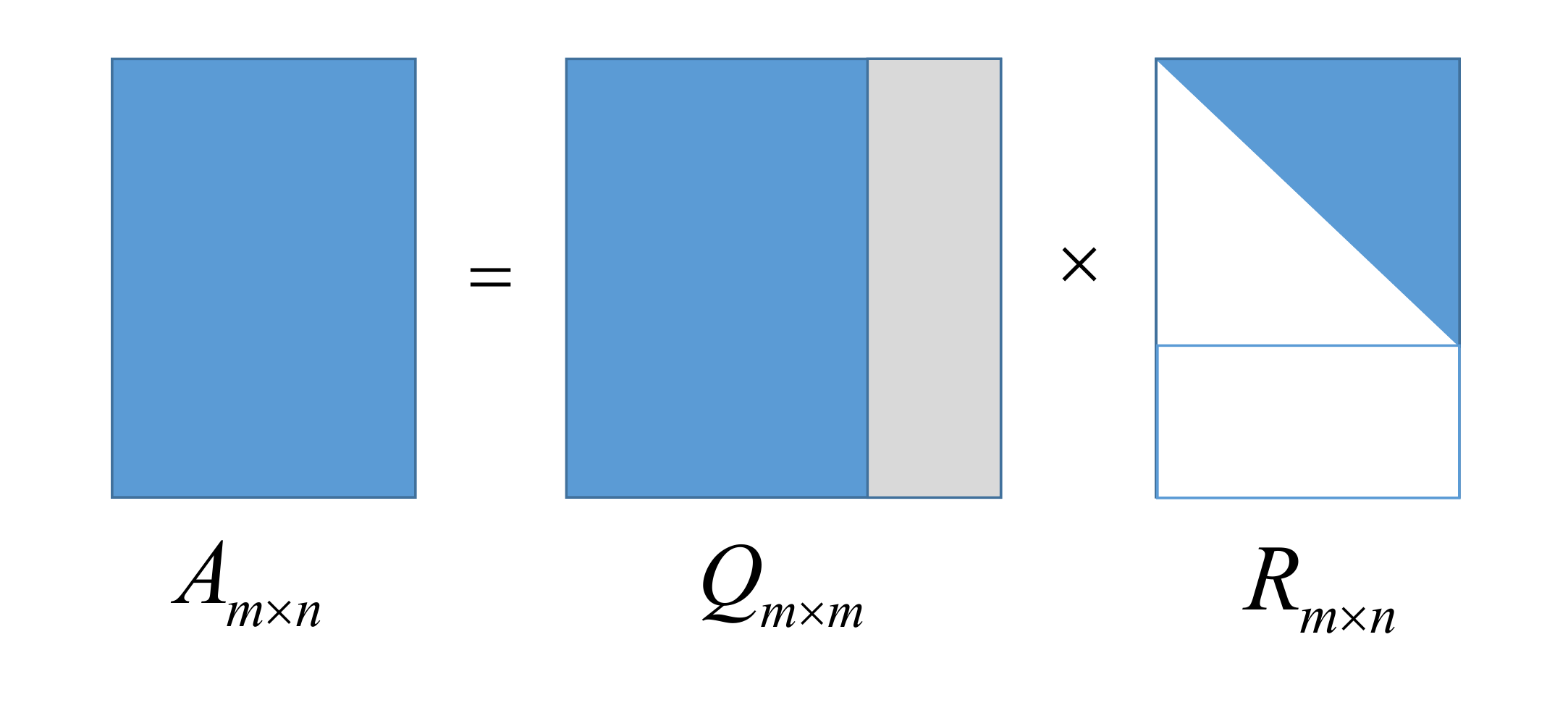}}
	\caption{Comparison of reduced and full QR decomposition.}
	\label{fig:qr-comparison}
\end{figure}

A full QR decomposition of an $m\times n$ matrix with independent columns goes further by appending additional $m-n$ orthonormal columns to $\bQ$ so that it becomes an $m\times m$ orthogonal matrix. In addition, rows of zeros are appended to $\bR$ so that it becomes an $m\times n$ upper triangular matrix. We call the additional columns in $\bQ$ as \textbf{silent columns} and additional rows in $\bR$ as \textbf{silent rows}. The comparison of the reduced QR decomposition and the full QR decomposition is shown in Figure~\ref{fig:qr-comparison} where silent columns in $\bQ$ are denoted in grey, blank entries are zero and blue entries are elements that are not necessarily zero.

\subsubsection{LQ Decomposition}
We previously proved the existence of the QR decomposition via the Gram-Schmidt process in which case we are interested in the column space of a matrix $\bA=[\ba_1, \ba_2, ..., \ba_n] \in \real^{m\times n}$. 
However, in many applications, we are also interested in the row space of a matrix $\bB=[\bb_1^\top; \bb_2^\top; ...;\bb_m^\top] \in \real^{m\times n}$, where $\bb_i$ is the $i$-th row of $\bB$ \citep{tanaka2003stochastic, schilders2009solution}. The successive spaces spanned by the rows $\bb_1, \bb_2, \ldots$ of $\bB$ are
$$
\cspace([\bb_1])\,\,\,\, \subseteq\,\,\,\, \cspace([\bb_1, \bb_2]) \,\,\,\,\subseteq\,\,\,\, \cspace([\bb_1, \bb_2, \bb_3])\,\,\,\, \subseteq\,\,\,\, \ldots,
$$
By applying QR decomposition on $\bB^\top = \bQ_0\bR$, we recover the LQ decomposition of the matrix $\bB = \bL \bQ$ where $\bQ = \bQ_0^\top$ and $\bL = \bR^\top$.

%
%
%

\subsubsection{UTV Decomposition: ULV and URV Decomposition}\label{section:ulv-urv-decomposition}

Finally, we come to the proof of the ULV decomposition. 
The existence of the ULV decomposition is from both the QR and LQ decomposition.
\begin{proof}[of Theorem~\ref{theorem:ulv-decomposition}]
	For any rank $r$ matrix $\bA=[\ba_1, \ba_2, \ldots, \ba_n]$, we can use a column permutation matrix $\bP$ such that the independent columns of $\bA$ appear in the first $r$ columns of $\bA\bP$. Without loss of generality, we assume $\bb_1, \bb_2, \ldots, \bb_r$ are the $r$ independent columns of $\bA$ and 
	$$
	\bA\bP = [\bb_1, \bb_2, \ldots, \bb_r, \bb_{r+1}, \ldots, \bb_n].
	$$
	Let $\bZ = [\bb_1, \bb_2, \ldots, \bb_r] \in \real^{m\times r}$. Since any $\bb_i$ is in the column space of $\bZ$, we can find a matrix $\bE\in \real^{r\times (n-r)}$ such that 
	$$
	[\bb_{r+1}, \bb_{r+2}, \ldots, \bb_n] = \bZ \bE.
	$$
	That is, 
	$$
	\bA\bP = [\bb_1, \bb_2, \ldots, \bb_r, \bb_{r+1}, \ldots, \bb_n] = \bZ
	\begin{bmatrix}
		\bI_r & \bE
	\end{bmatrix},
	$$
	where $\bI_r$ is an $r\times r$ identity matrix. Moreover, $\bZ\in \real^{m\times r}$ has full column rank such that its full QR decomposition is given by $\bZ = \bU\begin{bmatrix}
		\bR \\
		\bzero
	\end{bmatrix} \in \real^{m\times r}$, where $\bR\in \real^{r\times r}$ is an upper triangular matrix with full column rank and full row rank (because of the special structure of $\bR$) and $\bU$ is an orthogonal matrix. This implies 
	\begin{equation}\label{equation:ulv-smpl}
		\bA\bP = \bZ
		\begin{bmatrix}
			\bI_r & \bE
		\end{bmatrix}
		=
		\bU\begin{bmatrix}
			\bR \\
			\bzero
		\end{bmatrix}
		\begin{bmatrix}
			\bI_r & \bE
		\end{bmatrix}
		=
		\bU\begin{bmatrix}
			\bR & \bR\bE \\
			\bzero & \bzero 
		\end{bmatrix}.
	\end{equation}
	Since $\bR$ has full row rank, this means 
	$\begin{bmatrix}
		\bR & \bR\bE 
	\end{bmatrix}$ also has full row rank such that its full LQ decomposition is given by 
	$\begin{bmatrix}
		\bL & \bzero 
	\end{bmatrix} \bV_0$ where $\bL\in \real^{r\times r}$ is a lower triangular matrix and $\bV_0$ is an orthogonal matrix. Substitute into Equation~\eqref{equation:ulv-smpl}, we have 
	$$
	\bA = \bU \begin{bmatrix}
		\bL & \bzero \\
		\bzero & \bzero 
	\end{bmatrix}\bV_0 \bP^{-1}.
	$$
	Let $\bV =\bV_0 \bP^{-1}$ which is a product of two orthogonal matrices, and therefore is also an orthogonal matrix. This completes the proof.
\end{proof}

\begin{figure}[H]
	\centering  
	\vspace{-0.35cm} 
	\subfigtopskip=2pt 
	\subfigbottomskip=2pt 
	\subfigcapskip=-5pt 
	\subfigure[Reduced ULV decomposition]{\label{fig:ulvhalf}
		\includegraphics[width=0.47\linewidth]{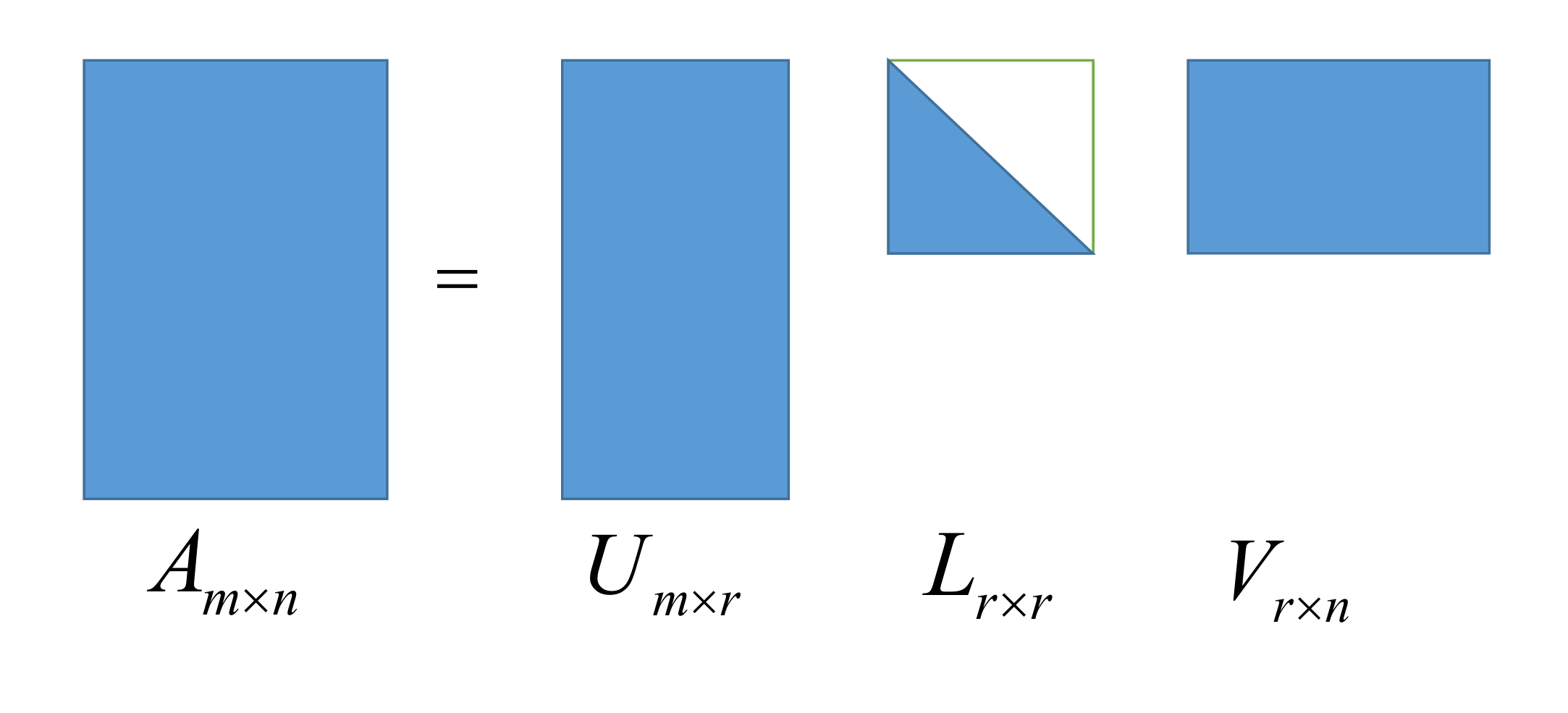}}
	\quad 
	\subfigure[Full ULV decomposition]{\label{fig:ulvall}
		\includegraphics[width=0.47\linewidth]{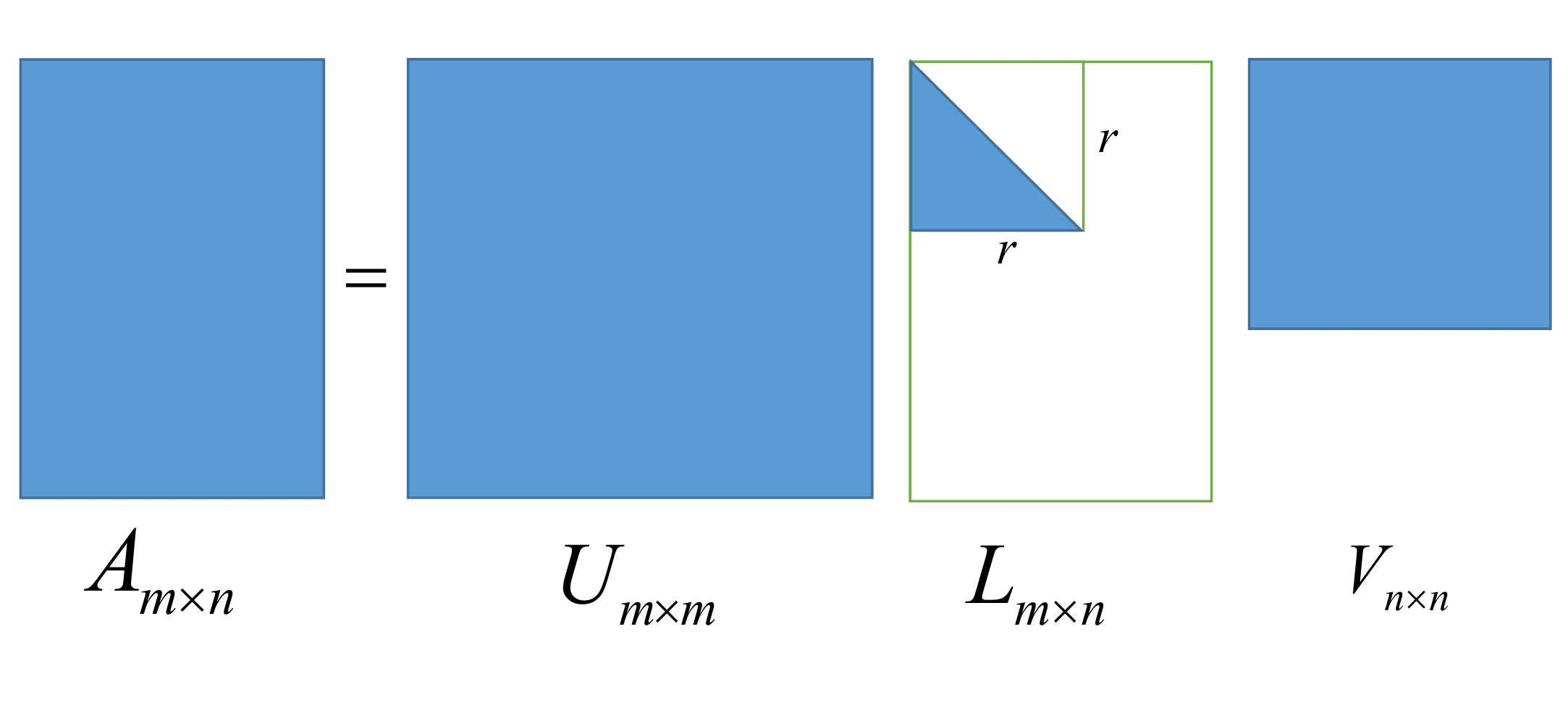}}
	\caption{Comparison of reduced and full ULV}
	\label{fig:ulv-comparison}
\end{figure}

Suppose the ULV decomposition of matrix $\bA$ is 
$$
\bA = \bU \begin{bmatrix}
	\bL & \bzero \\
	\bzero & \bzero 
\end{bmatrix}\bV.
$$
Let $\bU_0 = \bU_{:,1:r}$ and $\bV_0 = \bV_{1:r,:}$, i.e., $\bU_0$ contains only the first $r$ columns of $\bU$, and $\bV_0$ contains only the first $r$ rows of $\bV$. Then, we still have $\bA = \bU_0 \bL\bV_0$. This is known as the \textbf{reduced ULV decomposition}.
The comparison between the reduced and the full ULV decomposition is shown in Figure~\ref{fig:ulv-comparison} where white entries are zero and blue entries are not necessarily zero.

Similarly, we can also claim the URV decomposition as follows.
\begin{mdframed}[hidealllines=\mdframehideline,backgroundcolor=\mdframecolor]
	\begin{theorem}[Full URV Decomposition]\label{theorem:urv-decomposition}
		Every $m\times n$ matrix $\bA$ with rank $r$ can be factored as 
		$$
		\bA = \bU \begin{bmatrix}
			\bR & \bzero \\
			\bzero & \bzero 
		\end{bmatrix}\bV,
		$$
		where $\bU\in \real^{m\times m}$ and $\bV\in \real^{n\times n}$ are two orthogonal matrices, and $\bR\in \real^{r\times r}$ is an upper triangular matrix.
	\end{theorem}
\end{mdframed}
	The proof is similar to that of ULV, except this time, we need a row permutation rather than a column permutation. A detailed discussion can be found in \citep{lu2021numerical}.

Again, there is a version of reduced URV decomposition and the difference between the full and reduced URV can be implied from the context as shown in Figure~\ref{fig:ulv-comparison}. The ULV and URV sometimes are referred to as the UTV decomposition framework \citep{fierro1997low, golub2013matrix}.



\section{Gaussian Elimination Way}

CR decomposition is recently proposed in \citep{strang2021every, stranglu} that is a direct product of Gaussian elimination where the \textit{elementary transformations} will not change the rank of matrices. We firstly give the result and we will discuss the existence and the origin of this decomposition in the following paragraphs.
\begin{mdframed}[hidealllines=\mdframehideline,backgroundcolor=\mdframecolor]
	\begin{theorem}[CR Decomposition]\label{theorem:cr-decomposition}
		Any rank-$r$ matrix $\bA \in \real^{m \times n}$ can be factored as 
		$$
		\bA = \bC \bR
		$$
		where $\bC$ is some $r$ independent columns of $\bA$, and $\bR$ is a $r\times n$ matrix to reconstruct the columns of $\bA$ from columns of $\bC$. In particular, $\bR$ is the row reduced echelon form (RREF) of $\bA$ without the zero rows.
		
		The storage for the decomposition is then reduced or potentially increased from $mn$ to $r(m+n)$.
	\end{theorem}
\end{mdframed}

In short, we first compute the row reduced echelon form of matrix $\bA$ by $rref(\bA)$, Then $\bC$ is obtained by removing from $\bA$ all non-pivot columns (which can be determined by looking for columns in $rref(\bA)$ which do not contain a pivot). And $\bR$ is obtained by eliminating zero rows of $rref(\bA)$. And this is actually a special case of \textbf{rank decomposition} of matrix $\bA$. However, CR decomposition is so special that it involves in the row reduced echelon form. 

$\bR$ has a remarkable form whose $r$ columns containing the pivots form an $r\times r$ identity matrix. Note again that we can just remove the zero rows from the row reduced echelon form to obtain this matrix $\bR$. In \citep{strang2021every}, the authors give a specific notation for the row reduced echelon form without removing the zero rows as $\bR_0$:
$$
\bR_0 = rref(\bA)=
\begin{bmatrix}
	\bR \\
	\bzero
\end{bmatrix}=
\begin{bmatrix}
	\bI_r & \bF \\
	\bzero & \bzero
\end{bmatrix}\bP, \footnote{Permutation matrix $\bP$ on the right side of a matrix is to permute the columns of that matrix. 
}
$$
where the $n\times n$ permutation matrix $\bP$ puts the columns of $r\times r$ identity matrix $\bI_r$ into the correct positions, matching the first $r$ independent columns of the original matrix $\bA$.

The CR decomposition also reveals the great theorem of linear algebra that the row rank equals the column rank of any matrix. Thus we call it a \textit{Gaussian elimination} way.

\begin{proof}[of Theorem~\ref{lemma:equal-dimension-rank}, \textbf{The Third Way}]
	For CR decomposition of matrix $\bA=\bC\bR$, we have $\bR = [\bI_r, \bF ]\bP$, where $\bP$ is an $n\times n$ permutation to put the columns of the $r\times r$ identity matrix $\bI_r$ into the correct positions as shown above. It can be easily verified that the $r$ rows of $\bR$ are independent from the submatrix of $\bI_r$ (since $\bI_r$ is nonsingular) such that the row rank of $\bR$ is $r$. 
	
	Firstly, from the definition of CR decomposition, the $r$ columns of $\bC$ are from $r$ independent columns of $\bA$, the column rank of $\bA$ is $r$. Further, 
	
	$\bullet$ Since $\bA=\bC\bR$, all rows of $\bA$ are combinations of the rows of $\bR$. That is, the row rank of $\bA$ is no larger than the row rank of $\bR$;
	
	$\bullet$ From $\bA=\bC\bR$, we also have $(\bC^\top\bC)^{-1}\bC^\top\bC\bR = (\bC^\top\bC)^{-1}\bC^\top\bA$, that is $\bR = (\bC^\top\bC)^{-1}\bC^\top\bA$. $\bC^\top\bC$ is nonsingular since $\bC$ has full column rank $r$. Then all rows of $\bR$ are also combinations of the rows of $\bA$. That is, the row rank of $\bR$ is no larger than the row rank of $\bA$;
	
	$\bullet$ By ``sandwiching", the row rank of $\bA$ is equal to the row rank of $\bR$ which is $r$.
	
	Therefore, both the row rank and column rank of of $\bA$ are equal to $r$ from which the result follows.
\end{proof}

%
%
%
%

\subsection{Rank Decomposition}
We previously mentioned that the CR decomposition is a special case of rank decomposition. Formally, we prove the existence of the rank decomposition rigorously in the following theorem.

\begin{mdframed}[hidealllines=\mdframehideline,backgroundcolor=\mdframecolor]
	\begin{theorem}[Rank Decomposition]\label{theorem:rank-decomposition}
		Any rank-$r$ matrix $\bA \in \real^{m \times n}$ can be factored as 
		$$
		\bA = \bD \bF
		$$
		where $\bD \in \real^{m\times r}$ has rank $r$, and $\bF \in \real^{r\times n}$ also has rank $r$, i.e., $\bD,\bF$ have full rank $r$.
		
		The storage for the decomposition is then reduced or potentially increased from $mn$ to $r(m+n)$.
	\end{theorem}
\end{mdframed}
\begin{proof}[of Theorem~\ref{theorem:rank-decomposition}]
	By ULV decomposition in Theorem~\ref{theorem:ulv-decomposition}, we can decompose $\bA$ by 
	$$
	\bA = \bU \begin{bmatrix}
		\bL & \bzero \\
		\bzero & \bzero 
	\end{bmatrix}\bV.
	$$
	Let $\bU_0 = \bU_{:,1:r}$ and $\bV_0 = \bV_{1:r,:}$, i.e., $\bU_0$ contains only the first $r$ columns of $\bU$, and $\bV_0$ contains only the first $r$ rows of $\bV$. Then, we still have $\bA = \bU_0 \bL\bV_0$ where $\bU_0 \in \real^{m\times r}$ and $\bV_0\in \real^{r\times n}$. This is also known as the reduced ULV decomposition as shown in Figure~\ref{fig:ulv-comparison}. Let \{$\bD = \bU_0\bL$ and $\bF =\bV_0$\}, or \{$\bD = \bU_0$ and $\bF =\bL\bV_0$\}, we find such rank decomposition.
\end{proof}
The rank decomposition is not unique. Even by elementary transformations, we have 
$$
\bA = 
\bE_1
\begin{bmatrix}
	\bZ & \bzero \\
	\bzero & \bzero 
\end{bmatrix}
\bE_2,
$$
where $\bE_1 \in \real^{m\times m}, \bE_2\in \real^{n\times n}$ represent elementary row and column operations, $\bZ\in \real^{r\times r}$. The transformation is rather general, and there are dozens of these $\bE_1,\bE_2,\bZ$ matrices. Similar construction on this decomposition as shown in the above proof, we can recover another rank decomposition. 

Analogously, we can find such $\bD,\bF$ by SVD, URV, CR, Skeleton and many other factorization algorithms. However, we may connect the different rank decompositions by the following lemma.
\begin{mdframed}[hidealllines=\mdframehideline,backgroundcolor=\mdframecolor]
	\begin{lemma}[Connection Between Rank Decompositions]\label{lemma:connection-rank-decom}
		For any two rank decompositions of $\bA=\bD_1\bF_1=\bD_2\bF_2$, there exists a nonsingular matrix $\bP$ such that 
		$$
		\bD_1 = \bD_2\bP
		\qquad
		\text{and}
		\qquad 
		\bF_1 = \bP^{-1}\bF_2.
		$$
	\end{lemma}
\end{mdframed}
\begin{proof}[of Lemma~\ref{lemma:connection-rank-decom}]
	Since $\bD_1\bF_1=\bD_2\bF_2$, we have $\bD_1\bF_1\bF_1^\top=\bD_2\bF_2\bF_1^\top$. It is trivial that $rank(\bF_1\bF_1^\top)=rank(\bF_1)=r$ such that $\bF_1\bF_1^\top$ is a square matrix with full rank and thus is nonsingular. This implies $\bD_1=\bD_2\bF_2\bF_1^\top(\bF_1\bF_1^\top)^{-1}$. Let $\bP=\bF_2\bF_1^\top(\bF_1\bF_1^\top)^{-1}$, we have $\bD_1=\bD_2\bP$ and $\bF_1 = \bP^{-1}\bF_2$. 
\end{proof}

\subsection{Beyond UTV Decomposition}
We have shown that the UTV decomposition can be utilized to prove that the row and column ranks of a matrix are equal. From a basic course of linear algebra, the Gaussian elimination can results in the decomposition
$$
\bA = \bE_1 \begin{bmatrix}
	\bI_r & \bzero \\
	\bzero & \bzero 
\end{bmatrix}\bE_2,
$$
where $\bE_1, \bE_2$ are elementary transformations. By same argument as in Section~\ref{section:gaussian-way-main}, we can also find the same proof for Theorem~\ref{lemma:equal-dimension-rank}. 

Alternatively, we recall that the elementary row or column operations do not
alter the row rank or the column rank of $\bA$. Using these elementary transformations $\bE_1, \bE_2$, we will find the row and column ranks of $\bA$ are equal to those of $\begin{bmatrix}
	\bI_r & \bzero \\
	\bzero & \bzero 
\end{bmatrix}$ which proves the result as well.

\section{Skeleton/CUR Decomposition}
In Theorem~\ref{lemma:equal-dimension-rank}, we proved the row rank and the column rank of a matrix are equal. 
An elementary decomposition that is highly related to this theorem is called the skeleton decomposition.

\begin{mdframed}[hidealllines=\mdframehideline,backgroundcolor=\mdframecolor]
	\begin{theorem}[Skeleton Decomposition]\label{theorem:skeleton-decomposition}
		Any rank-$r$ matrix $\bA \in \real^{m \times n}$ can be factored as 
		$$
		\bA = \bC \bU^{-1} \bR,
		$$
		where $\bC$ is some $r$ independent columns of $\bA$, $\bR$ is some $r$ independent rows of $\bA$ and $\bU\in \real^{r\times r}$ is the nonsingular submatrix on the intersection. 
		
	\end{theorem}
\end{mdframed}

\begin{figure}[H]
	\centering
	\includegraphics[width=0.7\textwidth]{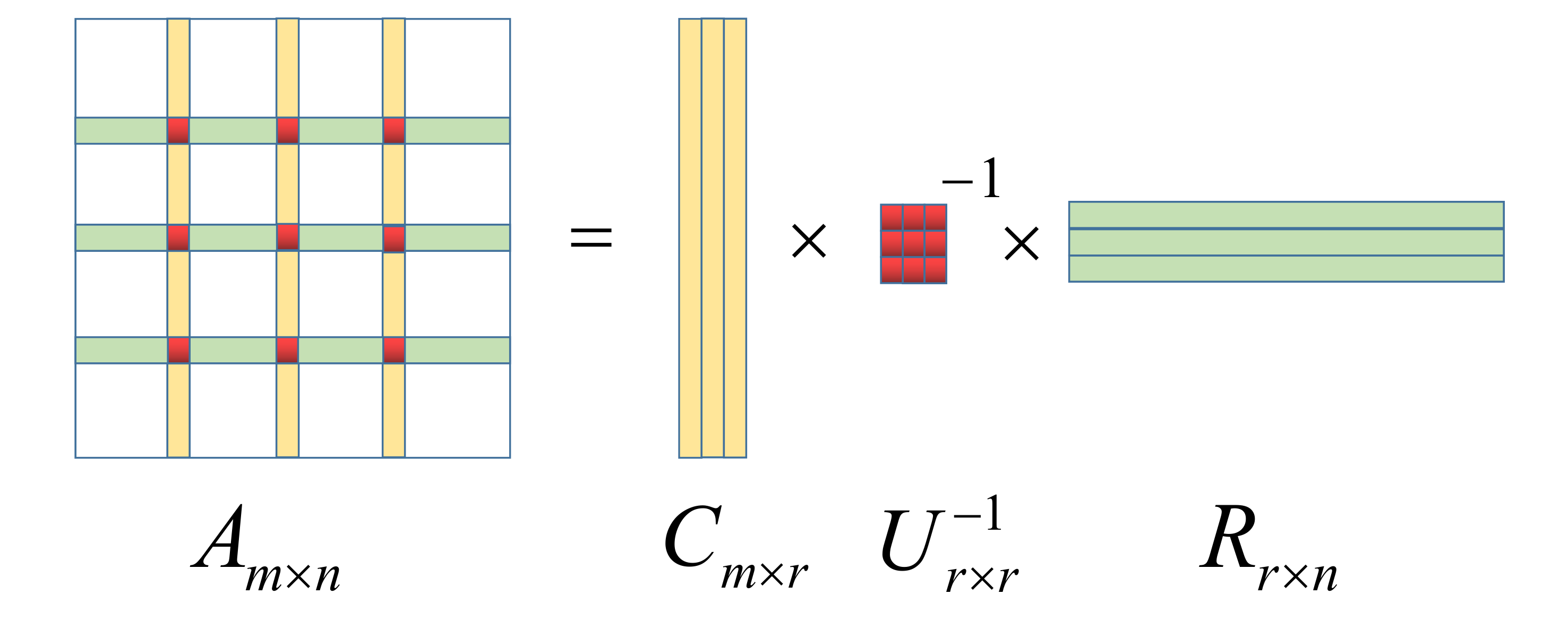}
	\caption{Demonstration of skeleton decomposition of a matrix}
	\label{fig:skeleton}
\end{figure}

Skeleton decomposition is also known as \textit{CUR decomposition}. The illustration of skeleton decomposition is shown in Figure~\ref{fig:skeleton}. Specifically, if $I,J$ are the indices of rows and columns selected from $\bA$ into $\bR$ and $\bC$ respectively, $\bU$ can be denoted as $\bU=\bA[I,J]$.

We have shown in previous sections that the Theorem~\ref{lemma:equal-dimension-rank} can be proved from Gram-Schmidt process or Gaussian elimination perspectives. The claim from Theorem~\ref{lemma:equal-dimension-rank} is also essential to prove the existence of skeleton decomposition. A moment of reflexion would reveal that Gram-Schmidt process or Gaussian elimination in turn can be utilized to produce such $\bC, \bU, \bR$ matrices for this decomposition (a detailed discussion is provided in \citep{lu2021numerical}).

We now prove the existence of the skeleton decomposition and the proof is rather elementary. 
\begin{proof}[of Theorem~\ref{theorem:skeleton-decomposition}]
	The proof relies on the existence of such nonsingular matrix $\bU$ which is central to this decomposition method. 
	
	\paragraph{Existence of such nonsingular matrix $\bU$} Since matrix $\bA$ is rank-$r$, we can pick $r$ columns from $\bA$ so that they are linearly independent. Suppose we put the specific $r$ independent columns $\ba_{i1}, \ba_{i2}, \ldots, \ba_{ir}$ into columns of an $m\times r$ matrix $\bN=[\ba_{i1}, \ba_{i2}, \ldots, \ba_{ir}] \in \real^{m\times r}$. The dimension of column space of $\bN$ is $r$ so that the dimension of row space of $\bN$ is also $r$ by Theorem~\ref{lemma:equal-dimension-rank}. Again, we can pick $r$ independent rows $\bn_{j1}^\top,\bn_{j2}^\top, \ldots, \bn_{jr}^\top $ from $\bN$ and put the specific $r$ rows into rows of a $r\times r$ matrix $\bU = [\bn_{j1}^\top; \bn_{j2}^\top; \ldots; \bn_{jr}^\top]\in \real^{r\times r}$. Using Theorem~\ref{lemma:equal-dimension-rank} again, the dimension of the column space of $\bU$ is also $r$ which means there are the $r$ independent columns from $\bU$. So $\bU$ is such a nonsingular matrix with size $r\times r$.
	
	As long as we find the nonsingular $r\times r$ matrix $\bU$ inside $\bA$, we can find the existence of the skeleton decomposition as follows.
	
	Since $\bU$ is a nonsingular matrix, the columns of $\bU$ are linearly independent. Thus the columns of matrix $\bC$ based on columns of $\bU$ are also linearly independent (i.e., select the $r$ columns of $\bA$ with same entries of the $\bU$ matrix. Here $\bC$ is equal to the $\bN$ we construct above). 
	
	As the rank of matrix $\bA$ is $r$, if we take any other column $\ba_i$ of $\bA$, $\ba_i$ can be represented as a linear combination of the columns of $\bC$, i.e., there exists a vector $\bx$ such that $\ba_i = \bC \bx$, for all $ i\in \{1, 2, \ldots, n\}$. Let $r$ rows of $\ba_i$ corresponding to the row entries of $\bU$ be $\br_i \in \real^r$ for all $i\in \{1, 2, \ldots, n\}$ (i.e., $\br_i$ contains $r$ entries of $\ba_i$). That is, select the $r$ entries of $\ba_i$'s corresponding to the entries of $\bU$ as follows:
	$$
	\bA = [\ba_1,\ba_2, \ldots, \ba_n]\in \real^{m\times n} \qquad \longrightarrow \qquad
	[\br_1, \br_2, \ldots, \br_n] \in \real^{r\times n}.
	$$
	Since $\ba_i = \bC\bx$, $\bU$ is a submatrix inside $\bC$, and $\br_i$ is a subvector inside $\ba_i$, we have $\br_i = \bU \bx$ which is equivalent to $\bx = \bU^{-1} \br_i$. Thus for every $i$, we have $\ba_i = \bC \bU^{-1} \br_i$. Combine the $n$ columns of such $\br_i$ into $\bR=[\br_1, \br_2, \ldots, \br_n]$, we obtain
	$$
	\bA = [\ba_1, \ba_2, \ldots, \ba_n] = \bC \bU^{-1} \bR,
	$$
	from which the result follows.
	
	In short, we first find $r$ independent columns of $\bA$ into $\bC\in \real^{m\times r}$. From $\bC$, we find a $r\times r$ nonsingular submatrix $\bU$. The $r$ rows of $\bA$ corresponding to entries of $\bU$ can help to reconstruct the columns of $\bA$. Again, the situation is shown in Figure~\ref{fig:skeleton}.
\end{proof}

In case $\bA$ is square and invertible, we have skeleton decomposition $\bA=\bC\bU^{-1} \bR$ where $\bC=\bR=\bU=\bA$ such that the decomposition reduces to $\bA = \bA\bA^{-1}\bA$.

\paragraph{CR decomposition vs skeleton decomposition} We note that CR decomposition and skeleton decomposition shares similar form. Even for the symbols used $\bA=\bC\bR$ for CR decomposition and $\bA=\bC\bU^{-1}\bR$ for skeleton decomposition. 

Both in CR decomposition and skeleton decomposition, we \textbf{can} select the first $r$ independent columns to obtain the matrix $\bC$ (the notation both in CR decomposition and in skeleton decomposition). So $\bC$ in CR decomposition and skeleton decomposition is exactly the same. On the contrary, $\bR$ in CR decomposition is the row reduced echelon form without the zero rows, and $\bR$ in skeleton decomposition is exactly some rows from $\bA$ so that $\bR$ has different meaning in the two decomposition methods. 

\paragraph{A word on the uniqueness of CR decomposition and skeleton decomposition} As mentioned above, both in CR decomposition and skeleton decomposition, we select the first $r$ independent columns to obtain the matrix $\bC$. In this sense, CR decomposition and skeleton decomposition has unique form. 
However, if we select the last $r$ independent columns, we will get a different CR decomposition or skeleton decomposition. We will not discuss this situation here as it is not the main interest of this text.  

To repeat, in the above proof for the existence of the skeleton decomposition, we first find the $r$ independent columns of $\bA$ into the $\bC$ matrix. From $\bC$, we find an $r\times r$ nonsingular submatrix $\bU$. From the submatrix $\bU$, we finally find the final row submatrix $\bR\in \real^{r\times n}$. A further question can be posed that if matrix $\bA$ has rank $r$ and matrix $\bC$ contains $r$ independent columns 
and matrix $\bR$ contains $r$ independent rows, then whether the $r\times r$ ``intersection" of $\bC$ and $\bR$ is invertible or not \footnote{We thank Gilbert Strang for raising this interesting question.}.

\begin{mdframed}[hidealllines=\mdframehideline,backgroundcolor=\mdframecolor]
	\begin{corollary}\label{corollary:invertible-intersection}
		If matrix $\bA \in \real^{m\times n}$ has rank $r$ and matrix $\bC$ contains $r$ independent columns 
		and matrix $\bR$ contains $r$ independent rows, then the $r\times r$ ``intersection" matrix $\bU$ of $\bC$ and $\bR$ is invertible.
	\end{corollary}
\end{mdframed}
\begin{proof}[of Corollary~\ref{corollary:invertible-intersection}]
	If $I,J$ are the indices of rows and columns selected from $\bA$ into $\bR$ and $\bC$ respectively, then, $\bR$ can be denoted as $\bR=\bA[I, :]$, $\bC$ can be represented as $\bC = \bA[:,J]$, and $\bU$ can be denoted as $\bU=\bA[I,J]$.
	
	Since $\bC$ contains $r$ independent columns of $\bA$, any column $\ba_i$ of $\bA$ can be represented as $\ba_i = \bC\bx_i = \bA[:,J]\bx_i$ for all $i \in \{1,2,\ldots, n\}$. This implies the $r$ entries of $\ba_i$ corresponding to the $I$ indices can be represented by the columns of $\bU$ such that $\ba_i[I] = \bU\bx_i \in \real^{r\times 1}$ for all $i \in \{1,2,\ldots, n\}$, i.e.,
	$$
	\ba_i = \bC\bx_i = \bA[:,J]\bx_i \in \real^{m} \qquad \longrightarrow  \qquad
	\ba_i[I] =\bA[I,J]\bx_i= \bU\bx_i \in \real^{r}.
	$$ 
	Since $\bR$ contains $r$ independent rows of $\bA$, the row rank and column rank of $\bR$ are equal to $r$. Combining the facts above, the $r$ columns of $\bR$ corresponding to indices $J$ (i.e., the $r$ columns of $\bU$) are linearly independent. 
	
	Again, by applying Theorem~\ref{lemma:equal-dimension-rank}, the dimension of the row space of $\bU$ is also equal to $r$ which means there are the $r$ independent rows from $\bU$, and $\bU$ is invertible. 
\end{proof}

\section{Acknowledgments}
We thank Gilbert Strang for raising the question formulated in the Corollary~\ref{corollary:invertible-intersection}, for a stream of ideas and references about the three factorizations from the steps of elimination, and the generous sharing of the manuscript of \citep{strang2021three}.


%

%
%



\vskip 0.2in
\bibliography{bib}

\end{document}